\documentclass[a4paper,12pt]{article}
\usepackage{amsmath,amsfonts,amsthm}
\usepackage{euscript}
\usepackage{mathbbol}
\usepackage{natbib}
\begin{document}
\title{Conciliation of Bayes and Pointwise \\
Quantum State Estimation:\\
asymptotic information bounds\\
in quantum statistics\footnote{This paper has appeared as R.D. Gill (2008), Conciliation of Bayes and pointwise
quantum state estimation
pp.\ 239--261 in {\it Quantum Stochastics and Information: Statistics, Filtering and Control}, V.P. Belavkin and M. Guta, eds., World Scientific. It was originally submitted perhaps in 2003 to {\it Annals of Statistics} under the title ``Asymptotic information bounds in quantum statistics''. It was accepted subject to minor corrections but with the suggestion to extend it with explanatory and background material. Unfortunately I let the dead-line pass. This version appeared in conference proceedings mentioned above. Finally revised and extended, co-author Madalin Guta, arXiv:1112.2078, published 2013 in another conference proceedings.}
}
\author{Richard D. Gill\thanks{URL: {\tt www.math.leidenuniv.nl/$\sim$gill}.
    Also affiliated with CWI, Amsterdam, the Netherlands, {\tt
      www.cwi.nl}.}  \\ Mathematical Institute, Leiden University \\ The
  Netherlands}
\date{19 Dec 2005 (arXiv v.1); this is arXiv version 5}

\newtheorem{theorem}{Theorem}
\newtheorem{question}{Question}
\newtheorem{fact}{Fact to remember}

\maketitle

\begin{abstract}
We derive an asymptotic lower bound on the Bayes risk when 
$N$ identical quantum systems whose state depends on a
vector of unknown parameters are jointly measured in an arbitrary way
and the parameters of interest estimated on the basis of 
the resulting data. The bound is an integrated version of a
quantum Cram\'er-Rao bound due to \citet{holevo82}, and it thereby links
the fixed $N$ exact  Bayesian optimality usually pursued in the physics literature
with the pointwise asymptotic optimality favoured in classical mathematical
statistics. By heuristic arguments the bound can be expected to be 
sharp. This does turn out to be the case in various important examples,
where it can be used to prove asymptotic optimality of 
interesting and useful measurement-and-estimation schemes. 
On the way we obtain a new family of
``dual Holevo bounds''  of independent interest.
\end{abstract}

\section{Introduction}

The aim of this paper is to derive asymptotic information bounds for
``quantum i.i.d.\ models''  in quantum statistics. That is to say, one has $N$ 
copies of a quantum system each in the same state depending on
an unknown vector of parameters $\theta$, and one wishes to estimate $\theta$,
or more generally a vector function of the parameters $\psi(\theta)$,
by making some measurement on the $N$ systems together. 
This yields data whose distribution depends on $\theta$ and on the choice of 
the measurement. Given the measurement, we therefore have a 
classical parametric statistical model, though not necessarily an i.i.d.\
model, since we are allowed to bring the $N$ systems together
before measuring the resulting joint system as one quantum object.
In that case the resulting data need not consist of (a function of)
$N$ i.i.d.\ observations, and a key quantum feature is that we
can generally extract more information about $\theta$ using
such ``collective'' or ``joint'' measurements than when we measure
the systems separately.
What is the best we can do as $N\to\infty$, when we are allowed to
optimize both over the measurement and over the ensuing
data-processing?

A heuristic, statistically motivated, approach to deriving methods
with good properties for large $N$ is to choose the measurement to 
optimize the Fisher information in the data, leaving it to the statistician to
process the data efficiently, using for instance maximum
likelihood or related methods, including Bayesian. This heuristic 
principle has already been shown to work in a number of special 
cases in quantum statistics. Since the measurement
maximizing the Fisher information typically
depends on the unknown parameter value this often has to
be implemented in a two-step approach, first using a small 
fraction of the $N$ systems to get a first approximation to
the true parameter, and then optimizing on the remaining
systems using this rough guess.

The approach favoured by many physicists is to choose a prior
distribution and loss function on grounds of symmetry and
physical interpretation, and then to \emph{exactly}
optimize the Bayes risk over all measurements and estimators, 
for any given $N$.
This approach succeeds in producing attractive methods
on those rare occasions 
when a felicitous combination of all the mathematical ingredients
leads to a simple and analytically tractable solution.
Now it has been observed in a number of problems that
the two approaches result in asymptotically equivalent
estimators, though the measurement schemes can be
strikingly different. Heuristically, this can be understood to follow
from the fact that, in the physicists' approach, for large $N$ 
the prior distribution should become increasingly irrelevant
and the Bayes optimal estimator close to the maximum
likelihood estimator.  Moreover, we expect those estimators
to be asymptotically normal with variances corresponding to
inverse Fisher information. 

Here we link the two approaches by
deriving a sharp asymptotic lower bound on the Bayes risk
of the physicists' approach, in terms of the optimal
Fisher information of the statisticians' approach.
This enables us to conclude the asymptotic optimality
of some heuristically motivated measurement-and-estimation 
schemes by showing that they attain the asymptotic bound.
Sometimes one can find in this way  asymptotically
optimal solutions which are much easier to implement
than the exactly optimal solution of the physicists' approach.
On the other hand, it also shows (if only heuristically) 
that the physicists' approach, when successful, leads to procedures
which are \emph{asymptotically} optimal for other prior distributions
than those used in the computation, also for loss functions 
only locally equivalent to their loss function of choice, and
also asymptotically optimal in a pointwise rather than
a Bayesian sense.

We derive our main result by combining an existing
quantum Cram\'er-Rao bound \citep{holevo82}
with the van Trees inequality, a Bayesian Cram\'er-Rao
bound from classical statistics (\citealp{vantrees68,
gilllevit95}). The former can be interpreted
as a bound on the Fisher information in an arbitrary
measurement on a quantum system, the latter is a bound
on the Bayes risk (for a quadratic loss function) in terms
of the Fisher information in the data. This means that
our result and its proof can be understood without any
familiarity with quantum statistics. Of course, to appreciate
the applications of the result, some further appreciation
of ``what is a quantum statistical model''  is needed.
The paper contains a brief summary of this; for more
information the reader is referred to the  papers
of \citet{barndorffnielsengilljupp03}, and \citet{gill01}.
For an overview of the ``state of the art'' in quantum asymptotic
statistics see \citet{hayashi05} which reprints papers of many authors 
together with introductions by the editor.

Let us develop enough notation to
state the main result of the paper and compare it with the
comparable result from classical statistics. Starting on
familar ground with the latter, suppose we
want to estimate a function $\psi(\theta)$ of a
parameter $\theta$, both represented by real column vectors
of possibly different dimension, based on $N$ i.i.d.\ observations from
a distribution with Fisher information matrix $I(\theta)$.
Let $\pi$ be a prior density on the parameter space and
let $\widetilde G(\theta)$ be a symmetric positive-definite matrix
defining a quadratic loss function 
$l(\widehat\psi^{(N)},\theta)= (\widehat\psi^{(N)}-\psi(\theta))^\top \widetilde G(\theta) (\widehat\psi^{(N)}-\psi(\theta))$. 
(Later we will use $G(\theta)$, without the tilde, in the special case when $\psi$ is $\theta$ itself).
Define the mean square error matrix 
$V^{(N)}(\theta)=\mathbb E_\theta (\widehat\psi^{(N)}-\psi(\theta)) (\widehat\psi^{(N)}-\psi(\theta))^\top$
so that the risk can be written $R^{(N)}(\theta)=\mathrm{trace}\,\widetilde G(\theta)V^{(N)}(\theta)$.
The Bayes risk is $R^{(N)}(\pi)=\mathbb E_\pi \mathrm{trace}\,\widetilde GV^{(N)}$.
Here, $\mathbb E_\theta$ denotes expectation over the data for given $\theta$,
$\mathbb E_\pi$ denotes averaging over $\theta$ with respect to the prior $\pi$.
The estimator $\widehat\psi^{(N)}$ is completely arbitrary.
We assume the prior density to be smooth, compactly supported 
and zero on the smooth boundary of its support. Furthermore a certain
quantity roughly interpreted as ``information in the prior'' 
must be finite.
Then it is very easy to show \citep{gilllevit95}, using the van Trees inequality,
 that under
minimal smoothness conditions on the statistical model,
\begin{equation}\label{mainresult1}
\mathop{\lim\inf}\limits_{N\to\infty}N R^{(N)}(\pi) ~\ge ~\mathbb E_\pi \mathrm{trace} \, G I^{-1}
\end{equation}
where $G=\psi' \widetilde G\psi'^\top$ and $\psi'$ is the matrix of
partial derivatives of elements of $\psi$ with respect to those of $\theta$.

Now in quantum statistics the data depends on the choice of measurement 
and the measurement should be tuned to the loss function.
Given a measurement $M^{(N)}$ on $N$ copies of
the quantum system, denote by $\overline I_M^{(N)}$ the average
Fisher information (i.e., Fisher information divided by $N$) in the
data. The \citet{holevo82} quantum Cram\'er-Rao bound, as 
extended by \citet{hayashimatsumoto04} to the quantum i.i.d.\ model, 
can be expressed as saying that, for all $\theta$, $G$, $N$ and  $M^{(N)}$,
\begin{equation}\label{holevo}
\mathrm{trace} \, G(\theta) (\overline I_M^{(N)}(\theta))^{-1}~\ge~\EuScript C_G(\theta)
\end{equation}
for a certain quantity $\EuScript C_G(\theta)$,
which depends on the specification of the quantum statistical model
(state of one copy,  derivatives of the state with respect to parameters, 
and loss function $G$)  \emph{at the point $\theta$} only, 
i.e., on local or pointwise model features (see (\ref{defCG1}) below).
According to as yet unpublished work of M.~Hayashi the bound is
asymptotically sharp. The idea behind his work is that locally, the quantum i.i.d.\ 
model is well approximated by a quantum Gaussian location model,
a quantum statistical problem for which 
the Holevo bound is sharp \citep{holevo82}.

We aim to prove that under minimal
smoothness conditions on the quantum statistical model, and conditions
on the prior similar to those needed in the classical case, 
but under essentially no conditions on the estimator-and-measurement sequence,
\begin{equation}\label{mainresult}
\mathop{\lim\inf}\limits_{N\to\infty}N R^{(N)}(\pi) ~\ge ~\mathbb E_\pi \EuScript C_G
\end{equation}
where, as before, $G=\psi' \widetilde G\psi'^\top$.
The main result (\ref{mainresult}) is exactly the bound one would hope for, 
from heuristic statistical principles, and one may also expect it
to be sharp, for the reasons mentioned above.
In specific models of interest, the right hand side is often easy
to calculate. Various specific measurement-and-estimator sequences,
motivated by a variety of approaches, can also be shown in
interesting examples to achieve the bound. The restrictions on
the prior can often be relaxed by approximating the prior of interest,
as we will show in our examples.

It was also shown in \citet{gilllevit95}, how---in the classical statistical context---one 
can replace a fixed prior $\pi$
by a sequence of priors indexed by $N$, concentrating more and more
on a fixed parameter value $\theta_0$, at rate $1/\sqrt N$. Following their
approach would, in the quantum context, lead to the pointwise asymptotic lower bounds 
\begin{equation}\label{regular}
\mathop{\lim\inf}\limits_{N\to\infty}N R^{(N)}(\theta) ~\ge ~\EuScript C_G(\theta)
\end{equation}
for each $\theta$, for \emph{regular}
estimators, and to local asymptotic minimax bounds 
\begin{equation}\label{minimax}
\lim_{M\to\infty}\mathop{\lim\inf}\limits_{N\to\infty}\sup_{\|\theta-\theta_0\|\le N^{-1/2}M}
N R^{(N)}(\theta) ~\ge ~\EuScript C_G(\theta_0)
\end{equation}
for \emph{all} estimators, but we do not further develop that
theory here.
In classical statistics the theory of Local Asymptotic Normality is the way
to unify, generalise, and understand this kind of result. We do not yet
have a theory of ``Q-LAN''  though there are indications that it may 
be possible to build such a theory. The results we obtain here
using more elementary tools do give  further support to the 
distant aim of building a Q-LAN theory.

The basic tools used  in this paper have now all been
mentioned, but as we shall see, the proof is not a routine
application of the van Trees inequality.
The missing ingredient will be provided by  the following new \emph{dual} 
bound to
(\ref{holevo}):
for all $\theta$, $K$, $N$ and  $M^{(N)}$,
\begin{equation}\label{dualholevo}
\mathrm{trace}\, K(\theta) \overline I_M^{(N)}(\theta)~\le~\EuScript C^K(\theta)
\end{equation}
where $\EuScript C^K(\theta)$ actually equals $\EuScript C_G(\theta)$ 
for a certain $G$ defined in terms of $K$ (as explained in
Theorem \ref{dualholevothm} below).
This is an \emph{upper} bound on Fisher information, in contrast to
(\ref{holevo}) which is a lower bound on inverse Fisher information.
The new inequality (\ref{dualholevo})
follows from the convexity of the sets of information matrices 
and of inverse information matrices  for arbitrary measurements 
on a quantum system, and these convexity
properties have a simple statistical explanation.
Such dual bounds have cropped
up incidentally in quantum statistics, for instance in \citet{gillmassar00},
but this is the first time a connection is established.

The argument for (\ref{dualholevo}), and given that, for (\ref{mainresult}), 
is based on some general structural features of
quantum statistics, and hence it is not necessary to be familiar with
the technical details of the set-up. 
In the next section we will summarize the
i.i.d.\ model in quantum statistics, focussing  on the
key facts which will be used in the proof of the dual Holevo bound 
(\ref{dualholevo}) and of our main result, the
asymptotic lower bound (\ref{mainresult}).
These proofs are given in a subsequent section, where
no further ``quantum''  arguments will be used. In a final section we will 
give three applications, leading to new results on some much
studied quantum statistical estimation problems.

\section{Quantum statistics: the i.i.d.\ parametric case.}

The basic objects in quantum statistics are \emph{states} and
\emph{measurements}, defined in terms of certain 
operators on a complex Hilbert space. 
To avoid technical complications we restrict attention to the 
finite-dimensional case, already rich in structure and
applications, when operators are represented by
ordinary (complex) matrices.

\paragraph{States and measurement}

The state of a $d$-dimensional system is represented by a $d\times d$
matrix $\rho$, called the \emph{density matrix} of the state, 
having the following properties: $\rho^*=\rho$ (self-adjoint
or Hermitian), $\rho\ge \mathbf 0$ (non-negative), $\mathrm{trace}(\rho)=1$
(normalized).  ``Non-negative'' actually implies ``self-adjoint'' but it does no harm to 
 emphasize both properties.  $\mathbf 0$ denotes the zero matrix; $\mathbf 1$ will denote 
 the identity matrix. 
 \medskip
 \\
 \noindent \emph{Example}: when $d=2$, every density matrix can be written
in the form $\rho=\frac12(\mathbf 1 + \theta_1 \sigma_1+\theta_2\sigma_2
+\theta_3 \sigma_3)$ where 
\begin{equation*}
\sigma_1=\left(\begin{matrix} 0& 1\\ 1& 0
\end{matrix}\right),\quad\sigma_2=\left(\begin{matrix} 0 &-i\\ i& 0
\end{matrix}\right),\quad\sigma_3=\left(\begin{matrix} 1 &0\\ 0 &-1
\end{matrix}\right)
\end{equation*}
 are the three Pauli matrices and where $\theta_1^2+\theta_2^2+\theta_3^2\le 1$. \hfill$\qed$
  \medskip
\\
  ``Quantum statistics'' concerns the situation when 
 the state of the system  $\rho(\theta)$ depends on a (column) vector $\theta$
 of $p$ unknown 
 (real) parameters.
  \medskip
 \\
 \emph{Example}: a completely unknown two-dimensional
 quantum state depends on a vector of three real parameters, 
 $\theta=(\theta_1,\theta_2,\theta_3)^\top$, known to lie in the unit ball.
 Various interesting submodels can be described geometrically: e.g.,
 the equatorial plane; the surface of the ball; a straight line through the origin.
 More generally, a completely unknown $d$-dimensional state depends on
 $p=d^2-1$ real parameters.\hfill $\qed$
  \medskip
 \\
\emph{Example}: in the previous example the two-parameter case 
obtained by demanding that  $\theta_1^2+\theta_2^2+\theta_3^2= 1$
is called the case of a two-dimensional pure state. In general,
a state is called pure if $\rho^2=\rho$ or equivalently $\rho$ has
rank one. A completely unknown pure $d$-dimensional state depends
on $p=2(d-1)$ real parameters.\hfill $\qed$
  \medskip
 
A measurement on a quantum system is characterized by the outcome space,
which is just a measurable space $(\EuScript X,\EuScript B)$, and a \emph{positive operator
valued measure} (POVM) $M$ on this space. This means that for each $B\in\EuScript B$  there corresponds a $d\times d$ non-negative self-adjoint matrix $M(B)$, 
together having the usual properties
 of an ordinary  (real) measure (sigma-additive), with moreover $M(\EuScript X)=\mathbf 1$.
 The probability distribution of the outcome of doing measurement $M$ on state $\rho(\theta)$
 is given by the Born law, or trace rule: $\Pr(\textrm{outcome}\in B)=\mathrm{trace}(\rho(\theta)M(B))$.
 It can be seen that this is indeed a bona-fide probability distribution on the sample
 space $(\EuScript X,\EuScript B)$. Moreover it has a density with respect to
 the finite real measure $\mathrm{trace}(M(B))$. 
 \medskip
 \\
\emph{Example}: the most simple measurement is defined by choosing an orthonormal
basis of $\mathbb C^d$, say $\psi_1$,\dots,$\psi_d$, taking the outcome space
to be the discrete space $\EuScript X=\{1,\dots,d\}$, and defining $M(\{x\})=\psi_x\psi_x^*$
for $x\in\EuScript X$;
or in physicists' notation, $M(\{x\})=|\psi_x\rangle\langle\psi_x|$.
One computes that $\Pr(\textrm{outcome}=x)=\psi_x^* \rho(\theta) \psi_x =
\langle \psi_x|\rho|\psi_x\rangle$. If the state is pure then $\rho=\phi\phi^*=
|\phi\rangle\langle\phi|$ for some $\phi=\phi(\theta)\in\mathbb C^d$ of length $1$ and depending
on the parameter $\theta$. One finds that  $\Pr(\textrm{outcome}=x)=|\psi_x^*\phi|^2
=|\langle\psi_x|\phi\rangle|^2$.   \hfill $\qed$
  \medskip

So far we have discussed state and measurement for a single quantum system.
This encompasses also the case of $N$ copies of the system,
via a tensor product construction, which we will now summarize.
The joint state of $N$ identical copies of a single system having state $\rho(\theta)$ is
$\rho(\theta)^{\otimes N}$, a density matrix on a space of dimension $d^N$. 
A joint or collective measurement on these systems is specified by a POVM on this
large  tensor product Hilbert space. An important point
is that joint measurements give many more possibilities than measuring the separate
systems independently, or even measuring the separate systems adaptively.

 \begin{fact} State plus measurement determines probability distribution of
 data.
 \end{fact}

\paragraph{Quantum Cram\'er-Rao bound.}

Our main input is going to be the \citet{holevo82} quantum Cram\'er-Rao bound, 
with its extension to the i.i.d.\ case due to \citet{hayashimatsumoto04}. 

Precisely because of quantum phenomena, different measurements, incompatible
with one another, are appropriate when we are interested in different components
of our parameter, or more generally, in different loss functions. The bound
concerns estimation of  $\theta$ itself rather than a function thereof, and depends on
a quadratic loss function defined by a symmetric real non-negative matrix $G(\theta)$ which 
may depend on the actual parameter value $\theta$. 
For a given estimator $\widehat\theta^{(N)}$ computed from the outcome
of some measurement $M^{(N)}$ on $N$ copies of our system, define its 
mean square error matrix $V^{(N)}(\theta)=\mathbb E_\theta (\widehat\theta^{(N)}-\theta)
(\widehat\theta^{(N)}-\theta)^\top$. The risk function when
using the quadratic loss determined by $G$ is 
$R^{(N)}(\theta)=\mathbb E_\theta (\widehat\theta^{(N)}-\theta)^\top G(\theta)  (\widehat\theta^{(N)}-\theta)=\mathrm{trace}(G(\theta)V^{(N)}(\theta))$. 

One may expect the risk of good measurements-and-estimators to decrease 
like $N^{-1}$ as $N\to \infty$.
The quantum Cram\'er-Rao bound confirms that this is the best
rate to hope for: it states that for unbiased estimators 
of a $p$-dimensional parameter $\theta$, based on
arbitrary joint measurements on $N$ copies, 
\begin{equation}\label{defCG1}
N R^{(N)}(\theta) ~\ge~  \EuScript C_G(\theta) ~=
~\inf_{\vec X,V:V \ge Z(\vec X)} \mathrm{trace}(G(\theta) V)
\end{equation}
where $\vec X = (X_1,\dots,X_p)$, the $X_i$ are
$d\times d$ self-adjoint matrices satisfying  $\partial/\partial\theta_i\, \mathrm{trace}(\rho(\theta) X_j)=\delta_{ij}$; 
$Z$ is the $p\times p$ self-adjoint matrix with elements $\mathrm{trace}(\rho(\theta)X_i X_j)$;
and $V$ is a real symmetric matrix.
It is possible to solve the optimization over $V$ for given $\vec X$ leading
to the formula
\begin{equation}\label{defCG2}
\EuScript C_G(\theta) ~=
~\inf_{\vec X} \mathrm{trace}\bigl( \Re (G^{1/2}Z(\vec X)G^{1/2} ) 
               + \mathrm{abs}\Im (G^{1/2}Z(\vec X )G^{1/2})\bigr)
 \end{equation}
where $G=G(\theta)$.
The absolute value of a matrix is found by diagonalising it and taking absolute values of
the eigenvalues. 
We'll assume that the bound is finite, i.e., there exists
$\vec X$ satisfying the constraints. A sufficient condition for this is that the
Helstrom quantum information matrix $H$ introduced in (\ref{helstrom})
below is nonsingular.

For specific interesting models, it often turns out not difficult to compute the bound
$\EuScript C_G(\theta)$. Note, it is a bound which depends only on the  
density matrix of one system ($N=1)$ and its derivative
with the respect to the parameter, and on the loss function, both at the given point $\theta$.
It can be found by solving a finite-dimensional optimization problem.

We will not be concerned with the specific form of the bound. What we are going to
need, are just two key properties. 

Firstly:  the bound is local, and applies to the larger
class of \emph{locally unbiased estimators}. This means to say that 
\emph{at the given point $\theta$}, $\mathbb E_\theta \widehat \theta^{(N)}=\theta$,
and at this point also $\partial/\partial\theta_i\, \mathbb E_\theta \widehat \theta_j^{(N)}=\delta_{ij}$.
Now, it is well known that the ``estimator'' $\theta_0+I(\theta_0)^{-1}S(\theta_0)$,
where $I(\theta)$ is Fisher information and $S(\theta)$ is score function, is locally
unbiased at $\theta=\theta_0$ \emph{and achieves the Cram\'er-Rao bound there}.
Thus the Cram\'er-Rao bound for \emph{locally} unbiased estimators is sharp.
Consequently, we can rewrite the  bound (\ref{defCG1}) in the form (\ref{holevo})
announced above,
where $\overline I_M^{(N)}(\theta)$ is the \emph{average} (divided by $N$) Fisher
information in the outcome of an arbitrary measurement $M=M^{(N)}$ on $N$ copies
and the right hand side is defined in (\ref{defCG1}) or (\ref{defCG2}).

\begin{fact} We have a family of computable lower bounds  
on the inverse average Fisher information matrix for an 
arbitrary measurement on $N$ copies, given by (\ref{holevo}) and (\ref{defCG1})
or (\ref{defCG2}),
\end{fact}

Secondly, for given $\theta$, define the following two sets of
positive-definite symmetric real matrices, in one-to-one correspondence
with one another through the mapping ``matrix inverse''. 
The matrices $G$ occurring in the definition are also 
taken to be positive-definite symmetric real.
\begin{equation}\label{defV}
\EuScript V = \{V:\mathrm{trace}(G V) \ge \EuScript  C_G~ \forall ~ G\},
\end{equation}
\begin{equation}\label{defI}
\EuScript I = \{I:\mathrm{trace}(G I^{-1}) \ge  \EuScript  C_G~ \forall ~ G\}.
\end{equation}
In the appendix to this paper, we give an algebraic proof
that that the set  $\EuScript I$ is convex (for $\EuScript V$, convexity is obvious),
and that the inequalities defining $\EuScript V$ define supporting hyperplanes 
to that convex set, i.e.,  all the inequalities are achievable in $\EuScript V$, 
or equivalently $\EuScript  C_G=\inf_{V\in\EuScript V} \mathrm{trace}(G V)$.

In fact, these properties have a statistical explanation, connected to
the fact that the quantum statistical problem of collective measurements
on $N$ identical quantum systems approaches a quantum Gaussian
problem as $N\to\infty$, see  \citet{gutakahn}.
It can be shown (\citealp{hayashi03}; Hayashi, personal communication; 
Gu\c t\u a, 2005, unpublished manuscript).
that \emph{$\EuScript V$ consists of all covariance 
matrices of locally unbiased estimators achievable (by suitable choice of measurement)
on a certain  $p$-parameter quantum Gaussian statistical model.
The inequalities defining $\EuScript V$ are the Holevo bounds
for that model, and each of those bounds is attainable}. Thus, for each $G$, 
there exists a $V\in\EuScript V$ achieving equality in 
$\mathrm{trace}(G V) \ge \EuScript  C_G$. 
It follows from this that $\EuScript I$ consists of all non-singular
information matrices
together with any non-singular matrix smaller than some information matrix, 
achievable by choice of measurement on the same quantum Gaussian model. 
Consider the set of information 
matrices attainable by some measurement together
with all smaller matrices; and consider the set of variance matrices of locally
unbiased estimators based on arbitrary measurements. Note that adding zero
mean noise to a locally unbiased estimator preserves its local unbiasedness,
so adding larger matrices to this set does not change it.
The set of information matrices is convex: choosing measurement $1$
with probability $p$ and measurement $2$ with probability $q$ (and
remembering your choice) gives a measurement whose Fisher information
is the convex combination of the informations of measurements $1$ and $2$.
Augmenting the set with all matrices smaller than something in the set, preserves
convexity. (The set of variances of locally unbiased estimators is convex, 
by a similar randomization argument). Putting this together, we obtain

\begin{fact}
For given $\theta$, 
both $\EuScript V$
and $\EuScript I$ defined in (\ref{defV}) and (\ref{defI}) are convex,
and all the inequalities defining these sets are achieved by points
in the sets.
\end{fact}

See the appendix for a direct algebraic proof.

\section{An asymptotic Bayesian information bound}

We will now introduce the van Trees inequality, a Bayesian Cram\'er-Rao
bound, and combine it with the Holevo bound (\ref{holevo}) via
derivation of a dual bound following from the convexity of the sets (\ref{defCG1}) 
and (\ref{defCG2}). We return to the problem of estimating the (real, column)  vector function 
$\psi(\theta)$ of the (real, column) vector parameter $\theta$
of a state $\rho(\theta)$ based on collective measurements of $N$ identical copies.
The dimensions of $\psi$ and of $\theta$ need not be the same. The sample
size $N$ is largely suppressed from the notation.
Let $V$ be the mean square error matrix of an arbitrary estimator $\widehat\psi$,
thus $V(\theta) = \mathbb E_\theta (\widehat\psi-\psi(\theta))(\widehat\psi-\psi(\theta))^\top$.
Often, but not necessarily, we'll have $\widehat\psi=\psi(\widehat\theta)$ for some
estimator of $\theta$.
Suppose we have a quadratic loss function 
$(\widehat\psi-\psi(\theta))^\top \widetilde G(\theta) (\widehat\psi-\psi(\theta))$
where $\widetilde G$ is a positive-definite matrix function of $\theta$, then
the Bayes risk with respect to a given prior $\pi$
can be written  $R(\pi)=\mathbb E_\pi \mathrm{trace}\, \widetilde G V$.
We are going to prove the following theorem:

\begin{theorem}
Suppose $\rho(\theta):\theta\in\Theta\subseteq \mathbb R^p$
is a smooth quantum statistical model and suppose
$\pi$ is a smooth prior density on a compact subset $\Theta_0\subseteq \Theta$, 
such that $\Theta_0$ has a piecewise smooth boundary, on which $\pi$ is zero.  
Suppose moreover the quantity $\EuScript J(\pi)$ defined in (\ref{defJpi}) below, is finite. Then
\begin{equation}\label{mainresult2}
\mathop{\lim\inf}\limits_{N\to\infty}N R^{(N)}(\pi) ~\ge ~\mathbb E_\pi \EuScript C_{G_0}
\end{equation}
where $G_0=\psi' \widetilde G\psi'^\top$ (and assumed to be positive-definite), 
$\psi'$ is the matrix of
partial derivatives of elements of $\psi$ with respect to those of $\theta$,
and $\EuScript C_{G_0}$  is defined by (\ref{defCG1})  or  (\ref{defCG2}).
\end{theorem}
``Once continuously differentiable'' is enough smoothness. Smoothness of
the quantum statistical model implies smoothness of the classical statistical
model following from applying an arbitrary measurement to $N$ copies of
the quantum state. Slightly weaker but more elaborate smoothness 
conditions on the statistical model and prior are spelled out
in \citet{gilllevit95}. The restriction that $G_0$ be non-singular
can probably be avoided by a more detailed analysis.

Let $\overline I_M$ denote the average Fisher information matrix
for $\theta$ based on a given collective measurement on the $N$ copies. Then the
 van Trees inequality states that for all matrix functions $C$ of $\theta$, of size
 $\mathrm {dim}(\psi)\times\mathrm{dim}(\theta)$, 
  \begin{equation}\label{e:vt}
 N\mathbb E_\pi  \mathrm{trace}\, \widetilde G V ~\ge~
 \frac{(\mathbb E_\pi\mathrm{trace}\, C\psi'^\top)^2}
 {\mathbb E_\pi \mathrm{trace}\, \widetilde G^{-1} C \overline I_M C^\top +\frac 1 N  
 \mathbb E_\pi \frac {(C\pi)'^\top \widetilde G^{-1} (C\pi)'} {\pi^2} }
 \end{equation}
 where the primes in $\psi'$ and in $(C\pi)'$ both denote differentiation,
 but in the first case converting the vector $\psi$ into the matrix of
 partial derivatives of elements of $\psi$ with respect to elements of
 $\theta$, of size $\mathrm {dim}(\psi)\times\mathrm{dim}(\theta)$,
 in the second case converting the matrix $C\pi$ into the
 column vector, of the same length as $\psi$,
 with row elements $\sum_j (\partial/\partial\theta_j)(C\pi)_{ij}$.
To get an optimal bound we need to choose $C(\theta)$ cleverly.
 
 First though,  note that the Fisher information appears in the denominator
 of the van Trees bound. This is a nuisance since we have a 
 Holevo's  \emph{lower}
 bound (\ref{holevo}) to the \emph{inverse} Fisher information.
 We would like to have an \emph{upper} bound on the
 information itself, say of the form (\ref{dualholevo}),
 together with a recipe for computing $\EuScript C^K$.
 
 All this can be obtained from the convexity of the sets
 $\EuScript I$ and $\EuScript V$ defined in (\ref{defI})
 and (\ref{defV}) and the non-redundancy of the
inequalities appearing in their definitions.
 Suppose $V_0$ is a boundary point of
$\EuScript V$. Define $I_0=V_0^{-1}$.  
Thus $I_0$ (though not necessarily an attainable average
information matrix $\overline I_M^{(N)}$) satisfies the Holevo
bound for each positive-definite $G$, and attains equality in one 
of them, say with $G=G_0$.  In the language of convex sets,
and ``in the $V$-picture'',
$\mathrm{trace}\, G_0 V= \EuScript C_{G_0}$ is a 
supporting hyperplane to $\EuScript V$ at $V=V_0$.

Under the mapping ``matrix-inverse" the hyperplane
$\mathrm{trace}\, G_0 V= \EuScript C_{G_0}$ in the $V$-picture 
maps to the smooth surface 
$\mathrm{trace}\, G_0 I^{-1}= \EuScript C_{G_0}$
touching the set $\EuScript I$ at $I_0$ in the $I$-picture.
Since $\EuScript I$ is convex, the tangent plane
to the smooth surface at $I=I_0$ must be 
a supporting hyperplane to $\EuScript I$
at this point.
The matrix derivative of the operation of matrix inversion
 can be written $\mathrm d A^{-1}/\mathrm d x= - A^{-1} (\mathrm d A / \mathrm d x) A^{-1}$.
This tells us that the equation of the tangent plane is
$\mathrm{trace}\, G_0 I_0 ^{-1} I I_0^{-1} = \mathrm{trace} \, G_0 I_0^{-1} =\EuScript C_{G_0}$.
Since this is simultaneously a supporting hyperplane to $\EuScript I$ we
deduce that for all $I\in \EuScript I$, 
$\mathrm{trace}\, G_0 I_0 ^{-1} I I_0^{-1} \le \EuScript C_{G_0}$. Defining
$K_0=I_0^{-1}G_0 I_0 ^{-1}$ and $\EuScript C^{K_0}=\EuScript C_{G_0}$ 
we rewrite this inequality as $\mathrm{trace}\, K_0 I \le \EuScript C^{K_0}$.

A similar story can be told when we start in the $I$-picture with a supporting
hyperplane (at $I=I_0$) to $\EuScript I$ of the form $\mathrm{trace}\, K_0 I = \EuScript C^{K_0}$
for some symmetric positive-definite $K_0$. It maps to the smooth surface
$\mathrm{trace}\, K_0 V^{-1} = \EuScript C^{K_0}$, with tangent plane
$\mathrm{trace}\, K_0 V_0^{-1} I V_0^{-1} = \EuScript C^{K_0}$ at $V=V_0=I_0^{-1}$.
By strict convexity of the function ``matrix inverse'', the tangent plane touches the
smooth surface only at the point $V_0$. Moreover, the smooth surface lies
above the tangent plane, but below $\EuScript V$. This makes $V_0$ the
unique minimizer of $\mathrm{trace}\, K_0 V_0^{-1} I V_0^{-1}$ in $\EuScript V$.

It would be useful to extend these computations to allow
singular $I$, $G$ and $K$.
Anyway, we summarize what we have so far in a theorem.
 \begin{theorem}\label{dualholevothm}
Dual to the Holevo family of lower bounds on average inverse
 information, $\mathrm{trace}\, G \overline I_M^{-1}\ge\EuScript C_G$ for each 
 positive-definite $G$, 
 we have a family of upper bounds on information, 
 \begin{equation}\label{dualholevo2}
 \mathrm{trace}\, K \overline I_M\le\EuScript C^K~~\text{for each}~~K.
 \end{equation}
 If $I_0\in\EuScript I$ satisfies $\mathrm{trace} \,G_0 I_0^{-1} = \EuScript C_{G_0}$
then with $K_0=I_0^{-1}G_0 I_0 ^{-1}$, $\EuScript C^{K_0}=\EuScript C_{G_0}$.
Conversely if $I_0\in\EuScript I$ satisfies $\mathrm{trace} \,K_0 I_0 = \EuScript C^{K_0}$
then with $G_0=I_0 K_0 I_0$, $\EuScript C_{G_0}=\EuScript C^{K_0}$. Moreover,
none of the bounds is redundant, in the sense that for all positive-definite
$G$ and $K$,
$\EuScript  C_G=\inf_{V\in\EuScript V} \mathrm{trace}(G V)$ and 
$\EuScript  C^K=\sup_{I\in\EuScript I} \mathrm{trace}(K I)$. 
The minimizer in the first equation is unique.
\end{theorem}

 Now we are ready to apply the van Trees inequality. First we make a guess for what the
 left hand side of (\ref{e:vt}) should look like, at its best. Suppose we use an estimator
 $\widehat\psi=\psi(\widehat\theta)$ where $\widehat\theta$ makes optimal
 use of the information in the measurement $M$. Denote
 now by $I_M$ the asymptotic normalized Fisher information of a sequence
 of measurements. Then we expect that the asymptotic normalized
 covariance matrix $V$
 of $\widehat\psi$ is equal to $\psi' I_M^{-1}\psi'^\top$ and therefore
 the asymptotic normalized Bayes risk should be 
 $\mathbb E_\pi\mathrm{trace}\,\widetilde G \psi' I_M^{-1}\psi'^\top
 =\mathbb E_\pi\mathrm{trace}\,\psi'^\top \widetilde G \psi' I_M^{-1}$. This is bounded
 below by the integrated Holevo bound $\mathbb E_\pi \EuScript C_{ G_0}$ with 
 $ G_0 = \psi'^\top \widetilde G \psi'$.
 Let $I_0\in\EuScript I$ satisfy
$\mathrm{trace}\,  G_0 I_0^{-1} =\EuScript C_{ G_0}$;
its existence and uniqueness are given by Theorem \ref{dualholevothm}.
(Heuristically we expect that $I_0$ is asymptotically attainable).
By the same Theorem, with $K_0=I_0^{-1}G_0 I_0^{-1}$,
 $\EuScript C^{K_0}=\EuScript C_{ G_0}=\mathrm{trace}\, G_0 I_0^{-1}
=\mathrm{trace}\,\psi'^\top \widetilde G \psi' I_0^{-1}$.

Though these calculations are informal, they lead us to try  the matrix function 
$C=\widetilde G \psi' I_0^{-1}$.  Define $V_0=I_0^{-1}$.
With this choice, in the numerator of the van Trees inequality,
we find the square of  $\mathrm{trace}\, C\psi'^\top=\mathrm{trace}\,\widetilde G \psi' I_0^{-1}\psi'^\top
=\mathrm{trace}\, G_0 V_0 =\EuScript C_{ G_0}$. In the main term of the denominator,
we find $\mathrm{trace}\,\widetilde G^{-1}\widetilde G \psi' I_0^{-1}\overline I_M  I_0^{-1} \psi'^\top 
\widetilde G=
\mathrm{trace}\,I_0^{-1} G_0 I_0^{-1} \overline I_M=
\mathrm{trace}\,K_0 \overline I_M \le\EuScript C^{K_0}=\EuScript C_{G_0}$
by the dual Holevo bound (\ref{dualholevo2}).
This makes the numerator of the van Trees bound equal to the square of
this part of the denominator, and using the inequality $a^2/(a+b)\ge a-b$ we find
  \begin{equation}\label{e:vt2}
 N\mathbb E_\pi  \mathrm{trace}\, G V ~\ge~
\mathbb E_\pi \EuScript C_{G_0} -
 \frac 1 N  
 \EuScript J(\pi)
  \end{equation}
where
\begin{equation}\label{defJpi}
\EuScript J(\pi) ~=~\mathbb E_\pi \frac {(C\pi)'^\top \widetilde G^{-1} (C\pi)'} {\pi^2}
\end{equation} 
with
$C=\widetilde G \psi' V_0$ and $V_0$ 
uniquely achieving in $\EuScript V$ 
the  bound $\mathrm{trace}\, G_0 V\ge\EuScript C_{G_0}$, where $ G_0 = \psi'^\top \widetilde G \psi'$.
Finally, provided $\EuScript J(\pi)$ is finite (which depends on the prior distribution
 and on properties of the model), we obtain the asymptotic lower bound
 \begin{equation}\label{e:vtasympt}
\mathop{\lim\inf}\limits_{N\to\infty} N\mathbb E_\pi  \mathrm{trace}\, \widetilde G V ~\ge~
\mathbb E_\pi\EuScript C_{G_0}. 
 \end{equation}

\section{Examples}


In the three examples discussed here, the loss function is
derived from a very popular (among the physicists) 
figure-of-merit in state estimation called
\emph{fidelity}.  Suppose we wish to estimate a state 
$\rho=\rho(\theta)$ by $\widehat\rho=\rho(\widehat\theta)$.
Fidelity measures the closeness of the
two states, being maximally equal to $1$ when the
estimate and truth coincide.
It is defined as $\mathrm{Fid}(\widehat\rho,\rho)=
\bigl(\mathrm{trace}(\sqrt {\rho^{\frac 12}\widehat \rho\rho^{\frac12}})\bigl)^2$
(some authors would call this \emph{squared} fidelity).
When both states are pure, thus $\rho=|\phi\rangle\langle\phi|$ and
$\widehat\rho=|\widehat\phi\rangle\langle\widehat\phi|$ where $\phi$ and $\widehat\phi$ are
unit vectors in $\mathbb C^d$, then $\mathrm{Fid}(\widehat\phi,\phi)
=|\langle \widehat\phi|\phi\rangle|^2$. There is an important characterization
of fidelity due to \citet{fuchs} which both explains its meaning and leads to many important
properties. Suppose $M$ is a measurement on the quantum system.
Denote by $M(\rho)$ the probability distribution of the outcome of
the measurement $M$ when applied to a state $\rho$. For two probability
distributions $P$, $\widehat P$ on the same sample space, let $p$ and
$\widehat p$ be their densities with respect to a dominating measure
$\mu$ and define the fidelity between these probability measures as
$\mathrm{Fid}(\widehat P,P)=\bigl(\int \widehat p^{\frac12} p^{\frac12}\mathrm d \mu\bigr)^2$.
In usual statistical language, this is the \emph{squared Hellinger affinity} between
the two probability measures. It turns out that
$\mathrm{Fid}(\widehat\rho,\rho)=\inf_M \mathrm{Fid}(M(\widehat \rho),M(\rho))$,
thus two states have small fidelity when there is a measurement which
distinguishes them well, in the sense that the Hellinger affinity between
the outcome distributions is small, or in other words, the $L_2$
distance between the root densities of the data under the two models
is large.  

Now suppose states are smoothly parametrized by a vector parameter $\theta$.
Consider the fidelity between two states with close-by parameter
values $\theta$ and $\widehat\theta$, 
and suppose they are measured with the same measurement $M$.
From the relation 
$ \int p^{\frac12}\widehat p^{\frac12}\mathrm d \mu = 1 -\frac 12 \| \widehat p^{\frac12} - p^{\frac12}\|^2$
and by a Taylor expansion to second order one finds
$1-\mathrm{Fid}(\widehat P,P)\approx\frac 14 (\widehat \theta -\theta)^\top
 I_M(\theta) (\widehat \theta -\theta)$ where $I_M(\theta)$ is the Fisher information in
 the outcome of the measurement $M$ on the state $\rho(\theta)$.
 We will define the \emph{Helstrom} quantum information matrix $H(\theta)$
 by the analogous relation 
\begin{equation}\label{helstrom}
1-\mathrm{Fid}(\widehat \rho,\rho)\approx\frac 14 (\widehat \theta -\theta)^\top
 H(\theta) (\widehat \theta -\theta).
 \end{equation}
It turns out that $H(\theta)$ is the smallest
 ``information matrix'' such that $I_M(\theta)\le H(\theta)$ for all measurements
 $M$.
 
Taking as loss function 
$l(\widehat\theta,\theta)=1-\mathrm{Fid}(\rho(\widehat \theta),\rho(\theta))$ we 
would expect (by a quadratic approximation to the loss) 
that $\mathbb E_\pi \EuScript C_{\frac14 H}$ is a
sharp asymptotic lower bound on $N$ times the Bayes risk. We will prove
this result for a number of special cases, in which by a fortuitous 
circumstance, the fidelity-loss function is \emph{exactly} quadratic
in a (sometimes rather strange) function of the parameter. 
The first two examples concern
a two-dimensional quantum system and are treated in depth in
\citet{baganetal05}; below we just outline some important
features of the application. In the second of those two examples
our asymptotic lower bound is an essential part of a proof of
asymptotic optimality of a certain measurement-and-estimation
scheme.

The third example concerns an unknown pure state of arbitrary
dimension. Here we are present a short and geometric 
proof of a surprising but little known  result of \citet{hayashi98} 
which shows that an extraordinarily simple measurement scheme
leads to  an asymptotically optimal estimator (providing the data
is processed efficiently). The analysis also
links the previously unconnected Holevo and Gill-Massar
bounds (\citealp{holevo82, gillmassar00}).

\subsection{Completely unknown spin half ($d$=2, $p$=3)}

Recall that a completely unknown $2$-dimensional quantum
state can be written 
$\rho(\theta)=\frac12(\mathbf 1 + \theta_1 \sigma_1+\theta_2\sigma_2
+\theta_3 \sigma_3)$, where $\theta$ lies in the unit ball in
$\mathbb R^3$. It turns out that  $\mathrm{Fid}(\widehat \rho,\rho)=
\frac 12 (1 + \widehat\theta \cdot \theta + (1-\|\widehat\theta\|^2)^{\frac12}(1-\|\theta\|^2)^{\frac12})$.
Define $\psi(\theta)$ to be the four-dimensional vector
obtained by adjoining $(1-\|\theta\|^2)^{\frac12}$ to $\theta_1$, $\theta_2$, $\theta_3$.
Note that this vector has constant length $1$.
It follows that  
$1-\mathrm{Fid}(\widehat \rho,\rho)=\frac14\|\widehat\psi-\psi\|^2$.
This is a quadratic loss-function for estimation of $\psi(\theta)$
with $\widetilde G = \mathbf 1$, the $4\times 4$ identity
matrix.
By Taylor expansion of both sides, we find that 
$\frac 14 H = \psi'^\top \widetilde G \psi' = G$
and conclude from Theorem 1 that $N$ times
$1-$ mean fidelity is indeed asymptotically lower bounded
by $\mathbb E_\pi \EuScript C_{\frac14 H}$.

In  \citet*{baganetal05} the exactly optimal measurement-and-estimation
scheme is derived and analysed in the case of a rotationally invariant
prior distribution over the unit ball. 
The optimal \emph{measurement} turns out not to
depend on the (arbitrary) radial part of the prior distribution,
and separates into two parts, one used for estimating 
the direction $\theta/\|\theta\|$, the other part for estimating
the length $\|\theta\|$. The Bayes optimal estimator of the length
of $\theta$ naturally depends on the prior. Because of these
simplifications it is feasible to compute the asymptotic value
of $N$ times the (optimal) Bayes mean fidelity, and this
value is $(3+2 \mathbb E_\pi\|\theta\|$)/4.

The Helstrom quantum information
matrix $H$ and the Holevo lower bound $\EuScript C_{\frac14 H}$ are also 
computed. 
It turns out that $\EuScript C_{\frac14 H}(\theta)=(3+2\|\theta\|)/4$.
Our asymptotic lower bound is not only correct but also, as expected,
sharp.

The van Trees approach does put some non-trivial conditions
on the prior density $\pi$. The most restrictive conditions are that
the density is zero at the boundary of its support and that the
quantity (\ref{defJpi}) be finite.
Within the unit ball everything is smooth, but there are
some singularities at the boundary of the ball. So our
main theorem does not apply directly to many priors
of interest. However there is an easy approximation
argument to extend its scope, as follows.

Suppose we start with a prior $\pi$ supported by the whole
unit ball which does not satisfy the conditions.  For any $\epsilon>0$
construct $\widetilde \pi=\widetilde \pi_\epsilon$ which is smaller than
$(1+\epsilon)\pi$ everywhere,  and $0$ for $\|\theta\|\ge 1-\delta$
for some $\delta>0$.
If the original prior $\pi$ is smooth enough we can arrange
that $\widetilde\pi$ satisfies the conditions of the van Trees
inequality, and makes (\ref{defJpi}) finite. $N$ times the Bayes risk
for $\widetilde\pi$ cannot exceed $1+\epsilon$ times 
that for $\pi$, and the same must also be true
for their limits. Finally, $\mathbb E_{\widetilde \pi_\epsilon}
\EuScript C_{\frac14 H}\to  \mathrm E_ \pi
\EuScript C_{\frac14 H}$ as $\epsilon\to 0$.

Some last remarks on this example: first of all, it is
known that \emph{only} collective measurements
can asymptotically achieve this bound. Separate measurements
on separate systems lead to strictly worse estimators.
In fact, by the same methods one can obtain the
sharp asymptotic lower bound $9/4$ (independent of the
prior), see Bagan, Ballester, Gill,
  Mu{\~n}oz-Tapia and Romero-Isart (2006b), when one allows
the measurement on the $n$th system to depend on
the data obtained from the earlier ones. Instead of
the Holevo bound itself, we use here a bound
of  \citet{gillmassar00}, which is actually has the form of
a dual Holevo bound.
(We give some more remarks on this at the end of the
discussion of the third example).
Secondly, our result gives strong heuristic support 
to the claim that the measurement-and-estimation
scheme developed in \citet*{baganetal05} 
for a specific prior and specific
loss function is also pointwise optimal in a minimax sense,
or among regular estimators, for loss functions which
are locally equivalent to fidelity-loss; and also asymptotically
optimal in the Bayes sense for other priors and locally
equivalent loss functions. In general, if the physicists' approach is
successful in the sense of generating a measurement-and-estimation
scheme which can be analytically studied and experimentally
implemented, then this scheme will have (for large $N$) good
properties independent of the prior and only dependent
on local properties of the loss.

\subsection{Spin half: equatorial plane ($d$=2, $p$=2)}

 \citet*{baganetal05}  also considered the case where
it is known that $\theta_3=0$, thus we now have a two-dimensional
parameter. The prior is again taken to be rotationally symmetric.
The exactly Bayes optimal measurement turns out (at least,
for some $N$ and for some priors)
to depend on the radial part of the prior. Analysis of
the exactly optimal measurement-and-estimation procedure
is not feasible since we do not know if this phenomenon persists
for all $N$. However there is a natural measurement,
which is exactly optimal for some $N$ and some priors, which
one might conjecture to be asymptotically optimal for all
priors. This sub-optimal measurement, combined with
the Bayes optimal estimator given the measurement,
can be analysed and it turns out that $N$ times
$1-$ mean fidelity converges to $1/2$ as $N\to \infty$,
independently of the prior. Again, 
the Helstrom quantum information
matrix $H$ and the Holevo lower bound $\EuScript C_{\frac14 H}$ are 
computed. 
It turns out that $\EuScript C_{\frac14 H}(\theta)=1/2$.
This time we can use our asymptotic lower bound to prove
that the natural sub-optimal measurement-and-estimator
is in fact asymptotically optimal for this problem.

For a $p$-parameter model the best one could every hope
for is that for large $N$ there are measurements with 
$\overline I_M$ approaching the Helstrom upper bound $H$.
Using this bound in the van Trees inequality gives the asymptotic
lower bound on $N$ times $1-$ mean fidelity of $p/4$.
The example here is a special case where this is attainable.
Such a model is called \emph{quasi-classical}.

If one restricts attention to separate measurements on
separate systems the sharp asymptotic lower bound is $1$,
twice as large, see Bagan, Ballester, Gill,
  Mu{\~n}oz-Tapia and Romero-Isart (2006b).

\subsection{Completely unknown $d$ dimensional pure state}

In this example we make use of the dual Holevo bound 
and symmetry arguments to show that in this example, the
original Holevo bound for a natural choice of $G$
(corresponding to fidelity-loss)
is attained by an extremely large class of measurements,
including one of the most basic measurements around,
known as ``standard tomography''.

For a pure state $\rho=|\phi\rangle\langle\phi|$, 
fidelity can be written $|\langle \widehat\phi|\phi\rangle|^2$
where $|\phi\rangle\in\mathbb C^d$ is a vector of unit length. 
The
state-vector can be multiplied by $e^{i a}$ for an arbitrary real phase
$a$ without changing the density matrix. The constraint of unit length
and the arbitrariness of the phase means that one can parametrize
the density matrix $\rho$ corresponding to $|\phi\rangle$ 
by $2(d-1)$ real parameters which we take to be our underlying
vector parameter $\theta$ (we have $d$ real parts and $d$ imaginary
parts of the elements of $|\phi\rangle$, but one constraint and one 
parameter which can be fixed arbitrarily).

For a pure state, $\rho^2=\rho$ so
$\mathrm{trace}(\rho^2)=1$.
Another way to write the fidelity in this case is as 
$\mathrm{trace}( \widehat\rho\rho)=\sum_{ij}(\Re(\widehat \rho_{ij})\Re(\rho_{ij})+
\Im(\widehat \rho_{ij})\Im(\rho_{ij}))$.  
So if we take $\psi(\theta)$ to be the vector
of length $2 d^2$ and of length $1$ containing the real and the imaginary parts of 
elements of $\rho$ we see that $1-\mathrm{Fid}(\widehat\rho,\rho)=\frac 12 \|\widehat\psi-\psi\|^2$.
It follows that $1-$ fidelity is a quadratic loss function in $\psi(\theta)$
with again $\widetilde G=\mathbf 1$. 

Define again the Helstrom quantum information matrix $H(\theta)$ for $\theta$
 by  $1-\mathrm{Fid}(\widehat \rho,\rho)\approx\frac 14 (\widehat \theta -\theta)^\top
 I_M(\theta) (\widehat \theta -\theta)$. Just as in the previous two examples
 we expect the asymptotic lower bound $\mathbb E_\pi \EuScript C_{\frac14 H}$
 to hold for $N$ times Bayes mean fidelity-loss, where $G=\frac14 H=
 \psi'^\top \widetilde G \psi' $.
 
 Some striking facts are known about estimation of a pure state. 
 First of all, from \citet{matsumoto02}, we know that the Holevo bound
 is attainable, for all $G$,  already at $N=1$. Secondly, from \citet{gillmassar00}
 we have the following inequality
 \begin{equation}\label{gillmassar}
 \mathrm{trace} H^{-1} \overline I_M\le d-1
 \end{equation}
 with \emph{equality} (in the case that the state is completely unknown) 
 for all \emph{exhaustive}
 measurements $M^{(N)}$ on $N$ copies of the state. Exhaustivity means,
 for a measurement with discrete outcome space, that $M^{(N)}(\{x\})$ 
 is a rank one matrix for
 each outcome $x$. The meaning of exhaustivity in general is by the same property for
 the density $m(x)$ of the matrix-valued measure $M^{(N)}$ with respect to a real
 dominating measure, e.g., $\mathrm{trace}(M^{(N)}(\cdot))$.
 This tells us that (\ref{gillmassar}) is one of the ``dual Holevo inequalities''.
 We can associate it with an original Holevo inequality once we know an 
 information matrix of a measurement attaining the bound. We will show that there
 is an information matrix of the form $\overline I_M=cH$ attaining the bound.
 Since the number of parameters (and dimension of $H$) is $2(d-1)$
 it follows by imposing equality in (\ref{gillmassar}) 
 that $c=\frac12$. The corresponding Holevo inequality 
 must be $ \mathrm{trace} \frac 12 H H^{-1} \frac 12 H \overline I_M^{-1}\ge d-1$
 which tells us that $\EuScript C_{\frac 14 H}=d-1$. 
 
 The proof uses an invariance property of the model. For any unitary matrix $U$
 (i.e., $UU^*=U^*U=\mathbf 1$) we can convert the pure state 
 $\rho$ into a new pure state $U\rho U^*$. The unitary matrices form a group
 under multiplication.  Consequently the group can be 
 thought to act on the parameter $\theta$ used to describe the pure state. Clearly the
 fidelity between two states (or the fidelity between their two parameters)
 is invariant when the same unitary acts on both states. This group action
 possesses the ``homogenous two point property'': for any two pairs of 
 states such that the fidelities between the members of each pair are the same, 
 there is a unitary transforming the first pair into the second pair.
 
 We illustrate this in the case $d=2$ where 
(first example, section 2), the pure
states can be represented by the surface of the unit ball in $\mathbb R^3$.
 It turns out that the action of the unitaries on the density matrices
 translates into the action of the group of orthogonal rotations on
 the unit sphere. Two points at equal distance on the sphere can
 be transformed by some rotation into any other two points at the
 same distance from one another; a constant distance between points
 on the sphere corresponds to a constant fidelity between the 
underlying states. 

In general, the pure states of dimension $d$ can be identified with
the Riemannian manifold $CP^{d-1}$ whose natural Riemannian
metric corresponds locally to fidelity (locally, $1-$ fidelity
is squared Riemannian distance) and whose isometries correspond
to the unitaries. This space posseses the homogenous two
point property, as we argued above. It is easy to show that the
\emph{only} Riemannian metrics invariant under isometries
on such a space are proportional to one another. Hence the
quadratic forms generating those metrics with respect to
a particular parametrization must also be proportional
to one another.

Consider a measurement whose outcome is actually an
estimate of the state, and suppose that this measurement
is \emph{covariant} under the unitaries. This means that
transforming the state by a unitary, doing the measurement
on the transformed state, and transforming the estimate back 
by the inverse of the same unitary, is the same 
(has the same POVM) as the original measurement.
The information matrix for such a measurement is generated
from the squared Hellinger affinity between the distributions
of the measurement outcomes under two nearby states,
just as the Helstrom information matrix is generated from
the fidelity between the states. If the measurement is
covariant then the Riemannian metric defined by the
information matrix of the measurement outcome must
be invariant under unitary transformations of the states.
Hence: \emph{the information matrix of any covariant measurement
is proportional to the Helstrom information matrix.}

Exhaustive covariant measurements certainly do exist. 
A particularly simple one is that, for each of the $N$ copies
of the quantum system, we independently and uniformly
choose a basis of $\mathbb C^d$ and perform the simple
measurement (given in an example in Section 2) 
corresponding to that basis.

The first conclusion of all this is: any exhaustive covariant
measurement has information matrix $\overline I_M^{(N)}$
equal to one half the Helstrom information matrix.
All such measurements attain the Holevo bound
$\mathrm{trace} \frac14 H ( \overline I_M^{(N)} )^{-1}\ge d-1$.
In particular, this holds for the i.i.d.\ measurement based
on repeatedly choosing a uniformly distributed random
basis of $\mathbb C^d$.

The second conclusion is that an asymptotic lower
bound on $N$ times $1-$ mean  fidelity is $d-1$. 
Now the exactly Bayes optimal measurement-and-estimation
strategy is known to achieve this bound. The measurement
involved is a mathematically elegant
collective measurement on the $N$ copies
together, but hard to realise in the laboratory. 
Our results show that one can expect to asymptotically
attain the bound by decent information processing
(maximum likelihood? optimal Bayes with uniform prior
and fidelity loss?)
following an arbitrary \emph{exhaustive covariant
measurement}, of which the most simple to implement
is the standard tomography measurement consisting of
an independent random choice of measurement basis
for each separate system.

In \citet{gillmassar00} the same bound as (\ref{gillmassar})
was shown to hold for separable (and in particular, for
adaptive sequential) measurements also in the mixed
state case. Moreover in the case $d=2$, any information
matrix satisfying the bound is attainable already at $N=1$.
This is used in \citet{baganetal06} to obtain sharp
asymptotic bounds to mean fidelity for separable
measurements on mixed qubits.

\section{Acknowledgements}

Thanks to Peter Jupp for advice on Riemannian manifolds, Manuel Ballester,
Madalin Gu\c{t}\u{a} and Jonas Kahn for their help in understanding the Holevo bound,
and to Emili Bagan, Ramon Mu\~nez-Tapia
and Alex Monras whose work led to the questions studied here.

I acknowledge financial support from the European Community project
 RESQ, contract no IST-2001-37559.
{


\bibliographystyle{Chicago}

\raggedright



}

\section*{Appendix: proof of convexity}

The first step is to show that 
\begin{equation}\label{defValt}
\EuScript V ~=~ \mathrm{clos}\{V: V\ge Z(\vec X)~\textrm{for some}~\vec X\}
\end{equation}
where, as before, $\vec X = (X_1,\dots,X_p)$, the $X_i$ are
$d\times d$ self-adjoint matrices satisfying  $\partial/\partial\theta_i\, \mathrm{trace}(\rho(\theta) X_j)=\delta_{ij}$; 
$Z$ is the $p\times p$ self-adjoint matrix with elements 
$\mathrm{trace}(\rho(\theta)X_i X_j)$;
and $V$ is a real symmetric matrix.

An easy computation
shows that $Z(p \vec X +(1-p) \vec Y)\le p Z( \vec X )+(1-p) Z( \vec Y)$
(check that the second derivative w.r.t.\ $p$ of $\langle \psi | Z(p \vec X +(1-p) \vec Y) | \psi\rangle$
is non-negative, for any complex vector $\psi$.)
This makes  $\{V: V\ge Z(\vec X)~\textrm{for some}~\vec X\}$, where $V$ is
self-adjoint, a convex set. Restricting to the real matrices in this set preserves
convexity, as does taking the closure of the set. By convexity, the definition (\ref{defCG1})
tells us that the equations $\mathrm{trace}(GV)=\EuScript C_G$ define supporting
hyperplanes to the set defined on the right hand side of (\ref{defValt}). 
Since a closed convex set is the intersection of the closed halfspaces defined
by its supporting hyperplanes, it follows that $\EuScript V$ as defined by
(\ref{defV}) can also be specified as (\ref{defValt}), and that all the
Holevo bounds $\mathrm{trace}(GV)\ge\EuScript C_G$ are attained
in $\EuScript V$. 

The convexity of $\EuScript I$, the set of inverses of elements of $\EuScript V$, 
is a lot more subtle. In the following argument I will
suppose that the state $\rho(\theta)$ is strictly positive. The proof is easily adapted
to the case of a model for a pure state. (More generally we need the notion of
D-invariant model and the $\EuScript L^2$ spaces defined by a quantum state, 
see \citealp{holevo82} or \citealp{hayashimatsumoto04}).

We can consider our model with $p$ parameters and a strictly positive
density matrix as a submodel of the model of
a completely unknown mixed state, which has $d^2-1$ parameters. 
Denote the parameter vector of the full model by $\phi$.
The submodel is parametrized by $\theta$, a subvector of
$\phi$. I'll use the terminology interest parameter, nuisance
parameter for the two subvectors of $\phi$ corresponding to
submodel parameters and auxiliary parameters.
Subscripts $1$, $2$ will be also used when we partition 
matrices or vectors according to these two parts.
By the strict positivity of $\rho$ we are working at a point in the interior
of the full model (this is one of the reasons why the argument 
needs to be adapted for a pure-state model).
Since $\mathrm{trace}\,\rho=1$, the partial derivatives of $\rho$ with respect to
the components of $\theta$ in submodel and $\phi$ in fullmodel
are traceless (i.e., have trace zero). It is easy to see from this 
that we may restrict the elements $X$ of $\vec X$, entering into
the Holevo bounds for the submodel,  and elements $Y$ of $\vec Y$,
entering into
the Holevo bounds for the full model,
to be such that $\mathrm{trace}\,\rho Y=0$. Such $Y$ form a $d^2-1$
dimensional real Hilbert space $\EuScript L^2_0(\rho)$ under the 
innerproduct $\langle X,Y\rangle_\rho = \Re\,\mathrm{trace}\,\rho X Y$.

Let $\rho'_i$ denote the partial derivative of $\rho$ with respect to
$\theta_i$ at the fixed parameter value under consideration. For
the submodel, define
the symmetric logarithmic derivatives $\lambda_i\in \EuScript L^2_0(\rho)$
by $\langle\lambda_i,X\rangle_\rho=\mathrm{trace}\,\rho'_i X$ for all
$X\in \EuScript L^2_0(\rho)$. The constraints $\mathrm{trace}\,\rho'_i X_j=\delta_{ij}$
translate into constraints $\langle\lambda_i,X_j\rangle_\rho=\delta_{ij}$ for all
$i,j\le p$.  
In the full model, I'll use the notation $\vec\mu$ for the vector of
symmetric logarithmic derivatives, and $\vec Y$ for a candidate vector
of $Y_i$, each of length $d^2-1$. Of course, $\vec\lambda$ is a subvector
of $\vec\mu$. In the full model, the constraints on $\vec Y$ translate into
$\langle\mu_i,Y_j\rangle_\rho=\delta_{ij}$ for all
$i,j\le d^2-1$.  The $\mu_i$ form a basis of $\EuScript L^2_0(\rho)$
of linearly independent  vectors.

Now in the full model,  the
constraints on the $Y_i$ make them uniquely defined. Thus for the full model,
the set $\EuScript V_{\mathrm{full}}$ is 
the set of all $(d^2-1)\times (d^2-1)$ real matrices  $W$
exceeding  the fixed self-adjoint matrix
$Z_{\mathrm{full}}=Z(\vec Y)$. Unfortunately, $Z_{\mathrm{full}}$
is singular. But we may describe $\EuScript I_{\mathrm{full}}$
as the closure of the set of all real matrices less than or equal
to $(Z_{\mathrm{full}}+\delta\mathbf 1)^{-1}$ for some $\delta>0$.
The convexity of both sets is trivial. This suggests that
we try to deal with the case of a $p$ parameter model by considering
it a submodel of the full $d^2-1$ parameter model.

The relation between inverse information matrices for full models and
submodels is complicated, but that between the information matrices
themselves is simple: the information matrix for a submodel is a
submatrix of the information matrix of a full model. Thus we might conjecture
that for every $I\in\EuScript I$, there exists a $W\ge 
Z_{\mathrm{full}}$ such that 
$I\le  (W^{-1})_{11}$, the subscript ``$11$''
indicating the submodel submatrix. However, it could be that
we have positive information for the submodel parameters,
but zero information for the auxiliary parameters. This would make
the corresponding inverse information matrix $W^{-1}$ for the full model
undefined. This problem can be solved by approximating
singular information matrices by nonsingular ones.
We will prove the following theorem:

\begin{theorem}
$V^{-1}\in\EuScript I$ if and only if there exist 
real matrices $W^{(n)}> Z_{\mathrm{full}}$,
with $((W^{(n)})^{-1})_{11}=(V^{(n)})^{-1}\to V^{-1}$ as $n\to\infty$.

In words, $\EuScript I$ is the closure of the set of 
$11$ submatrices of  real 
symmetric non-nonsingular matrices less
than or equal to  $(Z_{\mathrm{full}}+\delta\mathbf 1)^{-1}$ 
for some $\delta>0$. Consequently
$\EuScript I$ is convex.
\end{theorem}

\begin{proof}
The proof will work by frequent reparametrizations of the
nuisance part of the full model. By this we mean that
$\phi$ is transformed smoothly and one-to-one into,
say, $\psi$, in such a way that the interest component
of $\phi$ is unaltered. Under such a transformation,
the vector of symmetric logarithmic derivatives $\vec\mu$
transforms by premultiplication by an invertible matrix $C$
whose $11$ block is the identity and whose
$12$ block is zero, so the `interest'' part of $\vec\mu$ 
is unchanged. (Subject to $C$ being nonsingular, for which 
it is just necessary that the $22$
block is nonsingular, the $21$ block
of $C$ can be arbitrary). At the same time the vector
of operators $\vec Y$ transforms by premultiplication
by the transposed inverse of $C$. Consequently, 
$Z_{\mathrm{full}}$ is transformed into $(C^\top)^{-1}Z_{\mathrm{full}} C^{-1}$,
$W\ge Z_{\mathrm{full}}$ is transformed the same way, while
$W^{-1}$ is transformed into $C W^{-1} C^\top$.
We therefore see that the $11$ block (i.e., the submatrix
corresponding to 
the submodel) of  $W^{-1}$
\emph{remains invariant under reparametrization of the 
auxiliary or nuisance parameters}.

In the statement of the theorem the choice of parametrization
of the auxiliary parameters is arbitrary, and so can be chosen 
in any convenient way.
We take advantage of this possibility immediately,
in the proof of the the forwards implication of the theorem.

Suppose $V\ge Z(\vec X)$
for some $\vec X$ satisfying the usual constraints. Augment $\vec \lambda$ to
a vector $\vec \mu$ of $d^2-1$ linearly independent elements $\mu_i\in\EuScript L^2_0(\rho)$ 
such that $\langle \mu_i ,X_j\rangle_\rho=\delta_{ij}$ for all $i\le d^2-1$,
$j\le p$. (The extra elements can be an arbitrary basis of the orthocomplement of the $X_j$, 
it is easy to check that together with the old elements they are linearly independent,
hence because of their number, a basis). Next augment $\vec X$ to $\vec Y$, so that the
the orthogonality relation, with $X_j$ replaced by $Y_j$, also holds for $p<j\le d^2-1$.

For square matrices $A$, $B$ write $\mathrm{diag}(A,B)$ for the block
diagonal matrix with $A$ and $B$ as diagonal blocks corresponding to
interest and nuisance parts of the full model. 
Let $D_\epsilon=\mathrm{diag}(\mathbf 1,\epsilon \mathbf 1)$, this is the diagonal matrix 
with $1$'s on the interest parameter part of the
diagonal, $\epsilon$'s on the nuisance part.

We have $D_\epsilon Z_{\mathrm{full}} D_\epsilon \to \mathrm{diag}(Z(\vec X),\mathbf 0)\le 
\mathrm{diag}(V,\mathbf 0)$ as $\epsilon\to 0$. Therefore, for each $\epsilon>0$ we can find
$\delta=\delta(\epsilon)>0$ such that $D_\epsilon Z_{\mathrm{full}} D_\epsilon <
 \mathrm{diag}(V,\mathbf 0) +\delta \mathbf 1$ and moreover such that
 $\delta\to 0$ as $\epsilon\to 0$. Thus for each $\epsilon$,
 $Z_{\mathrm{full}} <
 D_\epsilon ^{-1}(\mathrm{diag}(V,\mathbf 0)+\delta \mathbf 1 )D_\epsilon^{-1}=W_\epsilon$
 where $((W_\epsilon)^{-1})_{11}\to V^{-1}$ as $\epsilon\to 0$.
 
 Choosing a sequence $\epsilon_n\to 0$ as $n\to \infty$ we have 
found $W^{(n)} > Z_{\mathrm{full}}$ for all $n$ with 
 $((W^{(n)})^{-1})_{11}\to V^{-1}$ as $n\to\infty$. Going back to the original parametrization
 does not alter $((W^{(n)})^{-1})_{11}$ so the forwards implication of the theorem is proved.
 
 Now for the backwards implication. Suppose I am given $W > Z_{\mathrm{full}}$,
 $(W^{-1})_{11}= V^{-1}$. Reparametrize
 the nuisance part of the full model so that $(W^{-1})_{12}=\mathbf 0$. This does not alter
 $(W^{-1})_{11}$ but does alter both interest and nuisance parts of $\vec Y$. 
 Denote the interest part of
 the transformed $\vec Y$ by $\vec X$. The inequality $W > Z_{\mathrm{full}}$
 remains true after the transformation, hence $W_{11} > Z(\vec X)$.
 Since $W$ is block diagonal, 
 we obtain from this
 $(W^{-1})_{11}\le (Z(\vec X)+\delta\mathbf 1)^{-1}$
 for some $\delta>0$. Taking the closure completes the proof.
 \end{proof}

 \end{document}